# Almost Locally Free Groups and the Genus Question

Anthony M. Gaglione
U.S. Naval Academy, Annapolis, MD 21402
Dennis Spellman
St. Joseph's University, Philadelphia, PA 19131

ABSTRACT. Sacerdote [Sa] has shown that the non-Abelian free groups satisfy precisely the same universal-existential sentences Th($F_2$)∩∀∃ in a first-order language $L_o$ appropriate for group theory. It is shown that in every model of Th($F_2$)∩∀∃ the maximal Abelian subgroups are elementarily equivalent to locally cyclic groups (necessarily nontrivial and torsion free). Two classes of groups are interpolated between the non-Abelian locally free groups and Remeslennikov's ∃-free groups. These classes are the **almost locally free groups** and the **quasi-locally free groups**. In particular, the almost locally free groups are the models of Th($F_2$)∩∀∃ while the quasi-locally free groups are the ∃-free groups with maximal Abelian subgroups elemenatarily equivalent to locally cyclic groups (necessarily nontrivial and torsion free). Two principal open questions at opposite ends of a spectrum are: (1.) Is every finitely generated almost locally free group free? (2.) Is every quasi-locally free group almost locally free? Examples abound of finitely generated quasi-locally free groups containing nontrivial torsion in their Abelianizations. The question of whether or not almost locally free groups have torsion free Abelianization is related to a bound in a free group on the number of factors needed to express certain elements of the derived group as a product of commutators.

## 1. Preliminaries

Let $L_o$ be the first-order language whose only relation symbol is = always interpreted as the identity relation and whose only function and constant symbols are a binary operation symbol $\cdot$ , a unary operation symbol $^{-1}$ and a constant symbol 1. $L_o$ shall be **the language of group theory**. Every sentence of $L_o$ is logically equivalent to one in **prenex normal form**

$$Q_1 x_1 \cdots Q_n x_n \psi(x_1, ..., x_n)$$

where for each i = 1,...,n, $Q_i$ is a quantifier (∀ or ∃) and $\psi(x_1, ..., x_n)$ is a formula of $L_o$ containing no quantifiers and containing free at most the distinct variables $x_1, ..., x_n$. Vacuous quantifications are permitted and it is agreed that each of $\forall x \psi$ and $\exists x \psi$ is logically equivalent to $\psi$ if the variable $x$ does not occur in $\psi$. A sentence of $L_o$ is both $\Pi_o$ and $\Sigma_o$ if it is logically equivalent to a sentence of $L_o$ containing no quantifiers.

A sentence of $L_o$ is $\Pi_1$ if it is logically equivalent to a **universal sentence**, i.e., a sentence of $L_o$ of the form $\forall x_1 \cdots \forall x_m \psi(x_1, ..., x_m)$. Here $\psi(x_1, ..., x_m)$ contains no quantifiers and contains at most the distinct variables $x_1, ..., x_m$. A sentence of $L_o$ is $\Sigma_1$ if it is logically equivalent to an **existential sentence**, $\exists x_1 \cdots \exists x_m \psi(x_1, ..., x_m)$ with similar provisos. Since the negation of a $\Pi_1$-sentence is a $\Sigma_1$-sentence and vice-versa, two groups G and H satisfy precisely the same $\Pi_1$-sentences of $L_o$ if and only if they satisfy precisely the same $\Sigma_1$-sentences of $L_o$. In that event, we say that G and H **have the same universal theory**.

A $\Pi_2$-sentence of $L_o$ is one logically equivalent to a **universal-existential sentence**, i.e., a sentence of the form $\forall x_1 \cdots \forall x_m \exists y_1 \cdots \exists y_n \psi(x_1, ..., x_m, y_1, \ldots, y_n)$. Here $\psi(x_1, ..., x_m, y_1, \ldots, y_n)$ contains no quantifiers and contains free at most the distinct variables $x_1, ..., x_m, y_1, \ldots, y_n$. A $\Sigma_2$-sentence of $L_o$ is one logically equivalent to an **existential-universal sentence**, i.e., a sentence of the form

$\exists x_1 \cdots \exists x_m \forall y_1 \cdots \forall y_n \psi(x_1, ..., x_m, y_1, \ldots, y_n)$ with similar provisos. Writing $\forall$, $\exists$, $\forall\exists$ and $\exists\forall$ for the set of $\Pi_1$-sentences, $\Sigma_1$-sentences, $\Pi_2$-sentences and $\Sigma_2$-sentences of $L_o$ respectively, we have since vacuous quantifications are permitted that $\forall \cup \exists \subseteq \forall\exists$.

Given a group G we write $Th(G)$ for the (first-order) **theory** of G (with respect to $L_o$). That is, $Th(G)$ is the set of all sentences of $L_o$ true in G. If $\alpha$ is a cardinal, then $F_\alpha$ shall be a group free of rank $\alpha$. Tarski has conjectured that if $\alpha$ and $\beta$ are any two cardinals each at least 2, then $Th(F_\alpha) = Th(F_\beta)$. Vaught [Gr, Theorem 4, p. 237] proved that Tarski's conjecture is correct if both $\alpha$ and $\beta$ are infinite. The case of finite rank remains open. Without taking a position on Tarski's conjecture, we let

$$\sum = \bigcap_{2 \leq r \leq \omega} Th(F_r)$$

where $\omega = \{0,1,2,\ldots\}$ is the first infinite cardinal. Thus, $\sum$ is the set of all sentences of $L_o$ true in every non-Abelian free group. Sacerdote [Sa] proved that every two non-Abelian free groups satisfy precisely the same $\Pi_2$- and $\Sigma_2$-sentences of $L_o$. Thus, $\Sigma \cap \forall\exists = Th(F_2) \cap \forall\exists$ is the set of all $\Pi_2$-sentences of $L_0$ true in every non-Abelian free group. Moreover since $\forall \cup \exists \subseteq \forall\exists$ and it is well known that the non-Abelian free groups do satisfy precisely the same universal and existential sentences of $L_o$, we must have

$$\Phi = \sum \cap (\forall \cup \exists) = Th(F_2) \cap (\forall \cup \exists)$$

contained in $Th(F_2) \cap \forall\exists$.

A model G of $Th(F_2) \cap \forall\exists$ satisfies all of the $\Pi_2$-sentences of $L_o$ true in every non-Abelian free group; however, there is no a priori prohibition against G satisfying a $\Pi_2$-sentence of $L_o$ false in some non-Abelian free group. By contrast, a model of $\Phi$ must satisfy precisely the same $\Pi_1$- and $\Sigma_1$-sentences of $L_o$ as the non-Abelian free groups. If S is a consistent set of sentences of $L_o$, then we write $\mathcal{M}(S)$ for the class of all models of S. We note that $\mathcal{M}(\Sigma) \subseteq \mathcal{M}(\text{Th}(F_2) \cap \forall\exists) \subseteq \mathcal{M}(\Phi)$. Remeslennikov calls a model of $\Phi$ a $\exists$-**free group** and we shall use "model of $\Phi$" and "$\exists$-free group" interchangeably.

**Definition 1.** *(B. Baumslag [B]): A group G is **fully residually free** or $\omega$-**residually free** provided to every finite nonempty subset $S \subseteq G - \{1\}$ there is a free group $F_S$ and an epimorphism $\varphi_S$: $G \to F_S$ such that $\varphi_S(g) \neq 1$ for all $g \in S$.*

Note that it follows from the Nielsen-Schreier Subgroup Theorem that full residual freeness is inherited by subgroups.

Our first proposition is a consequence of Theorem 2, p.279 of Grätzer [Gr].

**Proposition 2.** *Let S be a consistent set of sentences of $L_o$. Then there is a consistent set T of $\Pi_2$-sentences of $L_o$ such that $\mathcal{M}(S) = \mathcal{M}(T)$ if and only if $\mathcal{M}(S)$ is closed under the formation of direct unions.*

In particular, since $\Phi \subseteq \text{Th}(F_2) \cap \forall\exists$ each of $\mathcal{M}(\Phi)$ and $\mathcal{M}(\text{Th}(F_2) \cap \forall\exists)$ is closed under the formation of direct unions.

**Proposition 3.** *(Remeslennikov [Re1]): A finitely generated group is $\exists$-free if and only if it is non-Abelian and fully residually free.*

We easily deduce from Propositions 2 and 3 the following

**Corollary 4.** *A group is $\exists$-free if and only if it is non-Abelian and locally fully residually free.*

Remeslennikov has shown that if G is $\exists$-free then G is $\Lambda$-free (in the sense of Bass [Ba]) for a suitable ordered Abelian group $\Lambda = (A, \leq)$ depending upon G; furthermore, he has shown that if G has finite rank k, then A can be taken to be free Abelian of rank not exceeding 3k. Consequences of this result include for a $\exists$-free group G the following:

- G is torsion free.

- Every 2-generator subgroup is either free or Abelian.

- G is **commutative transitive** in the sense that centralizers of nontrivial elements are Abelian (and hence maximal Abelian),

- G is **CSA** in the sense that its maximal Abelian subgroups M satisfy the malnormality condition: $g^{-1}Mg \cap M \neq 1 \Rightarrow g \in M$.

- If G is finitely generated, then G satisfies the maximal condition for Abelian subgroups.

Note that in any torsion free, commutative transitive group the maximal Abelian subgroups are **isolated** in the sense that if $g^n \in M$ for an integer n > 0, then $g \in M$. In particular, this must be true of ∃-free groups.

**Proposition 5.** *(V. Huber Dyson [D]): Every finitely presented, residually finite group has solvable word problem.*

Since it is well known that free groups are residually finite, we have the following

**Corollary 6.** *Every finitely presented, ∃-free group has solvable word problem.*

Two groups G and H are **elementarily equivalent** provided Th(G) = Th(H). It was shown in Gaglione and Spellman [GS] that in every model of Σ every maximal Abelian subgroup is elementarily equivalent to the infinite cyclic group. A set of necessary and sufficient conditions for two Abelian groups G and H to be elementarily equivalent was found by Szmielew [Sz]. Since we will need them, we reproduce these conditions here.

If G is an Abelian group, m > 0 is an integer and $\{g_i : i \in I\} \subseteq G$, then the $g_i$ will be said to be **linearly independent modulo m** provided whenever a finite product $\prod_i g_i^{e_i}$ is the identity element 1 of G one must have $e_i \equiv 0 \pmod{m}$ for all i, and $\{g_i : i \in I\}$ is said to be **linearly independent modulo m in the strong sense** provided whenever a finite product $\prod_i g_i^{e_i}$ lies in $G^m$ one must have $e_i \equiv 0 \pmod{m}$ for all i. If G is an Abelian group, p a prime, and k a positive integer, then three ranks

$$\rho^{(i)}[p,k](G) \in \varpi \cup \{\infty\} \ (i = 1,2,3)$$

are defined by

- $\rho^{(1)}[p,k](G) = n$ if n is the maximal number of elements of order $p^k$ and linearly independent modulo $p^k$.

- $\rho^{(2)}[p,k](G) = n$ if n is the maximal number of elements of G linearly independent modulo $p^k$ in the strong sense.

- $\rho^{(3)}[p,k](G) = n$ if n is the maximal number of elements of order $p^k$ and linearly independent modulo $p^k$ in the strong sense.

Two Abelian groups G and H are elementarily equivalent if and only if they are either both of finite exponent or both of infinite exponent and $\rho^{(i)}[p,k](G) = \rho^{(i)}[p,k](H)$, i = 1,2,3, for all primes p and positive integers k.

Suppose $A \neq 1$ is a multiplicatively written torsion free, locally cyclic group. A complete set of isomorphism invariants for A is given in Kurosh (see p. 207-210, [K]). A is determined up to isomorphism by (the equivalence class of - under a suitable relation) an infinite sequence $(a_1, a_2, \ldots, a_n, \ldots)$ all of the terms of which are either nonnegative integers or the symbol $\infty$. Writing $p_n$ for the nth positive rational prime in natural order, the equation $x^{p_n} = c$ will have a solution in A for arbitrary $c \in A$ if and only if $a_n = \infty$.

It is an exercise in number theory to show that if S(A) is the set of all $p_n$ such that $a_n = \infty$, then A is elementarily equivalent to the additive group of the ring $\mathsf{Z}[S(A)^{-1}]$. (Here we are using the notation of [L] on p.60.) We shall abuse notation and denote this additive group by this symbol $\mathsf{Z}[S(A)^{-1}]$. Thus every nontrivial, torsion free locally cyclic group is elementarily equivalent to $\mathsf{Z}[S^{-1}]$ for a suitable set S (possibly empty) of primes.

## 2. ALMOST LOCALLY FREE GROUPS

**Definition 7.** *A model G of Th($F_2$) $\cap \forall \exists$ will be called* **almost locally free**.

Since $\Pi_2$-sentences are preserved under the formation of direct unions, every non-Abelian locally free group is almost locally free. If $^*F_2$ is an ultrapower of $F_2 = <a_1, a_2; >$ with respect to a non-principal ultrafilter on then first infinite cardinal $\omega$, then $^*F_2$ is a model of Th($F_2$) $\cap \forall \exists$ (even more so, of $\Sigma$) which is not locally free. (See [BeSl] for a discussion of ultrapowers.)

**Theorem 8.** *Let G be an almost locally free group and let M be a maximal Abelian subgroup of G. Then there is a set S(M) of primes (possibly empty) such that M is elementarily equivalent to $\mathsf{Z}[S(M)^{-1}]$.*

**Proof:** For each integer n > 1 let $\varphi_n$ be the following $\Pi_2$-sentence of $L_o$:

$$\forall x_o \forall x_1 \forall x_2 \exists y (((x_o \neq 1) \wedge (x_o x_1 = x_1 x_o) \wedge (x_o x_2 = x_2 x_o))$$

$$\rightarrow ((x_o y = y x_o) \wedge \bigvee_{\substack{0 \leq i,j < n \\ (i,j) \neq (0,0)}} (x_1^i x_2^j = y^n))) \tag{1}$$

asserting that modulo the n-th powers any two elements of a maximal Abelian subgroup are linearly dependent. (Here our notion of dependence coincides with the negation of linearly independent modulo n in the strong sense.) Since $\varphi_n$ is true in every non-Abelian free group, it must be true in G. In particular if p is a prime and M is maximal Abelian in G, then the vector space $M/M^p$ over the p-element field is such that any two vectors are linearly dependent. It follows that either $M = M^p$ or $M/M^p$ is cyclic of order p.

Let $S(M)$ be the set of primes p such that $M = M^p$. It will suffice to show that for arbitrary primes p not in $S(M)$ and arbitrary positive integers k, it is the case that $M/M^{p^k}$ is cyclic of order $p^k$ since that is true of $\mathsf{Z}[S(M)^{-1}]/p^k\mathsf{Z}[S(M)^{-1}]$ for such primes p. So suppose p is a prime and $c_o \neq 1$ is an element of M which is not a p-th power.

We first claim that the $p^k$ elements $1, c_o, c_o^2, \ldots, c_o^{p^k-1}$ are pairwise incongruent modulo $M^{p^k}$. Suppose to deduce a contradiction that there are integers r and s with $0 \leq r < s < p^k$ such that $c_o^r \equiv c_o^s \pmod{M^{p^k}}$. Let $p^m = \gcd(s-r, p^k)$. Observe that $0 \leq m < k$ and that $t = \frac{s-r}{p^m}$ is not divisible by p. Furthermore since in a $\exists$-free group roots are unique when they exist, $c_o^t \equiv 1 \pmod{M^{p^{k-m}}}$. But since t is prime to p, t has an inverse modulo $p^{k-m}$ and $c_o \equiv 1 \pmod{M^{p^{k-m}}}$. But this is impossible since $m < k$ and $c_o$ is not a p-th power.

To finish our proof that $M/M^{p^k}$ is cyclic of order $p^k$, it will therefore suffice to show that if $d \in M$ is arbitrary, then there must exist an index $i = i(d)$ with $0 \leq i < p^k$ such that $d \equiv c_o^i \pmod{M^{p^k}}$. For each prime p and positive integer k, let $\psi_{p,k}$ be the following $\Pi_2$-sentence of $L_o$:

$$\forall x_o \forall x_1 \exists y ((x_o \neq 1) \wedge (x_o x_1 = x_1 x_o))$$

$$\rightarrow ((x_o y = y x_o) \wedge ((x_o = y^p) \vee \bigvee_{i=0}^{p^k-1} (x_1 = x_o^i y^{p^k}))) \tag{2}$$

whose interpretation in any free group F is that for each nontrivial element $x_o$ either $x_o$ is a p-th power or $x_o$ generates $Z_F(x_o)$ modulo $Z_F(x_o)^{p^k}$. (Here $Z_G(x)$ means the centralizer of an element x in a group G.) This is true in any non-Abelian free group since if $x_o \neq 1$ is not a p-th power, then it generates $Z_F(x_o)$ modulo $Z_F(x_o)^{p^k}$. Thus $\psi_{p,k}$ is true in the almost locally free group G for every prime p and positive integer k. In particular, since $c_o \neq 1$ and is not a p-th power in G it must generate M modulo $M^{p^k}$.

Hence, $M/M^{p^k}$ is cyclic of order $p^k$ for every prime p not in $S(M)$ and every positive integer k. It follows from Szmielew's conditions stated earlier that M is elementarily equivalent to the locally cyclic group $\mathsf{Z}[S(M)^{-1}]$. ♠

Of course the free Abelian group of rank 2, for example, cannot be elementarily equivalent to any locally cyclic group. Thus the model $<a,b,c; bc=cb>$ of $\Phi$ containing a maximal Abelian subgroup $<b,c>$, free Abelian of rank 2, cannot be almost locally free. Hence, the class of almost locally free groups is a proper subclass of the class of ∃-free groups. Indeed since the only finitely generated Abelian groups elementarily equivalent to $\mathsf{Z}[S^{-1}]$ for a set S of primes are the infinite cyclic groups (forcing S to be empty) and since finitely generated ∃-free groups satisfy the maximal condition for Abelian subgroups, the maximal Abelian subgroups in every finitely generated almost locally free group must be infinite cyclic. Every finitely generated, almost locally free group is **2-free** in the sense that every 2-generator subgroup is free. This is so because in any ∃-free group every pair of non-commuting elements freely generates a subgroup. In fact, finitely generated almost locally free groups are 3-free since it was shown in Fine, Gaglione, Myasnikov, Rosenberger and Spellman [FGMRS] that every 2-free, residually free group is 3-free.

**Definition 9.** *A **quasi-locally free group** is a ∃-free group which satisfies for each integer n>1 the $\Pi_2$-sentence $\varphi_n$:*

$$\forall x_o \forall x_1 \forall x_2 \exists y(((x_o \neq 1) \wedge (x_o x_1 = x_1 x_o) \wedge (x_o x_2 = x_2 x_o))$$
$$\to ((x_o y = y x_o) \wedge \bigvee_{\substack{0 \leq i,j < n \\ (i,j) \neq (0,0)}} (x_1^i x_2^j = y^n))) \tag{3}$$

*and for each prime p and positive integer k the $\Pi_2$-sentence $\psi_{p,k}$ :*

$$\forall x_o \forall x_1 \exists y(((x_o \neq 1) \wedge (x_o x_1 = x_1 x_o)) \to$$
$$((x_o y = y x_o) \wedge ((x_o = y^p) \vee \bigvee_{i=0}^{p^k - 1} (x_1 = x_o^i y^{p^k}))))) \tag{4}$$

(For interpretations of (3) and (4), see the remarks after (1) and (2), respectively.)

Equivalently, a quasi-locally free group is a ∃-free group in which every maximal Abelian subgroup is elementarily equivalent to a locally cyclic group. A third formulation is that the ∃-free group is quasi-locally free provided to every maximal Abelian subgroup M ≤ G there is a set S(M) of primes (possibly empty) such that M is elementarily equivalent to $\mathsf{Z}[S(M)^{-1}]$.

Thus we have the chain of inclusions illustrated by the following Venn diagram:

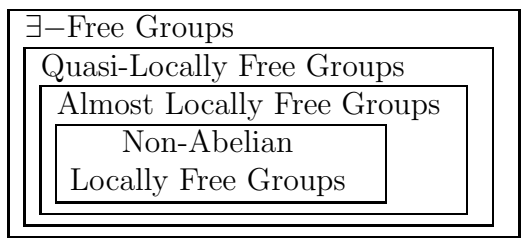

We do not know whether or not every quasi-locally free group is almost locally free. The non-orientable surface group

$$K = <a_1, a_2, a_3, a_4; a_1^2 a_2^2 a_3^2 a_4^2 = 1>$$

of genus 4 is an example of a quasi-locally free group since it is a finitely generated, non-Abelian, fully residually free group in which the maximal Abelian subgroups are infinite cyclic. Indeed to every finite set $g_1,\ldots, g_n$ of nontrivial elements of K there is a retraction $\varphi : K \to <a_1, a_2>$ such that $\varphi(g_i) \neq 1$ for all $i = 1,\ldots, n$. It follows that K and $F_2 = <a_1, a_2>$ satisfy precisely the same universal and existential sentences of the extended language $L_o[F_2]$ formed by adjoining the elements of $F_2$ as new constants. Thus, if $\exists \overline{x} \forall \overline{y} \psi(\overline{x}, \overline{y})$ is an existential-universal sentence of $L_o$, then for some tuple $\overline{c}$ of elements of $F_2$ the universal sentence $\forall \overline{y} \psi(\overline{c}, \overline{y})$ of $L_o[F_2]$ is true in $F_2$ so must also be true in K. It follows that $\exists \overline{x} \forall \overline{y} \psi(\overline{x}, \overline{y})$ is true in K whenever it is true in $F_2$. Put another way, K is a model of $Th(F_2) \cap \exists \forall$. However, we do not know whether or not K is a model of $Th(F_2) \cap \forall \exists$ nor even whether or not there are finitely generated, non-free almost locally free groups!

Thus,

Question 1: Is every finitely generated, almost locally free group free?

Question 2: Is every quasi-locally free group almost locally free?

We reiterate here questions posed in the same vein in [GS]:

Question 1′: Is every finitely generated model of $\Sigma$ free?

Question 2′: Is every ∃-free group in which the maximal Abelian subgroups are all elementarily equivalent to the infinite cyclic group a model of $\Sigma$?

We shall return in the next section to these questions. For now we further study the classes of almost locally free groups and quasi-locally free groups. We shall declare a subgroup H of the ∃-free group G to be centralizer closed in G provided for every $h \in H - \{1\}$ it is the case that

$Z_H(h) = Z_G(h)$. The the class of quasi-locally free groups is closed under the formation of non-Abelian, centralizer closed subgroups. This is so since if G is a non-Abelian, locally fully residually free group in which the centralizer of every nontrivial element is elementarily equivalent to a locally cyclic group, then so is every non-Abelian, centralizer closed subgroup. Furthermore, the free product of two quasi-locally free groups is again quasi-locally free. This is so since the class of ∃-free groups is closed under the formation of free products (see Theorem 3 [GS]); moreover, in a free product the maximal Abelian subgroups either lie in a conjugate of a free factor or are infinite cyclic. The free product of a quasi-locally free group with a group elementarily equivalent to Z $[S^{-1}]$ for some set S of primes or even of two groups elementarily equivalent respectively to Z $[S_1^{-1}]$ and Z $[S_2^{-1}]$ where $S_1$ and $S_2$ are sets of primes is again quasi-locally free. We cannot determine at this time whether or not the same assertions hold for the class of almost locally free groups. However, one closure property we can assert with confidence is common to both classes is closure under the formation of direct unions. This is so since the class of quasi-locally free groups, no less than the class of almost locally free groups, is a model class and has a set of $\Pi_2$-axioms.

To every set S of primes, we may construct the non-Abelian locally free group Z*Z$[S^{-1}]$. Thus a representative of every equivalence class of nontrivial torsion free locally cyclic groups is attainable as a maximal Abelian subgroup. However, we shall presently see that in every quasi-locally free group some maximal Abelian subgroup must be elementarily equivalent to the infinite cyclic group.

To that end, let X = $\{x_1, x_2\}$ and let w($x_1,x_2$) be a freely reduced, nonempty word on the alphabet X∪$X^{-1}$. We shall say that w($x_1,x_2$) is primitive provided it defines an element of a free basis for the free group on X.

**Proposition 10.** *Let G be a quasi-locally free group. Let g ∈ G-{1}. If there is a freely reduced, nonempty word w($x_1,x_2$) which is neither primitive nor a proper power and a pair of non-commuting elements a and b of G such that g = w(a,b), then the centralizer of g is elementarily equivalent to the infinite cyclic group.*

Proof: It suffices to show that g is not a p-th power in G for any prime p. Suppose to deduce a contradiction that g = $c^p$ for some c ∈ G and prime p. Consider the subgroup H = <a,b,c> of G. H is a homomorphic image of the one relator group

$$Y = <y_1, y_2, y_3; w(y_1,y_2) = y_3^p>.$$

H is non-Abelian since it contains the non-commuting elements a and b. Every non-Abelian subgroup of a ∃-free group is ∃-free; so, H, being ∃-free and finitely generated is fully residually free ( by Proposition 3). Also H is fully residually free of rank 2 by Lemma 8 [B]. But that contradicts a theorem of G. Baumslag and A. Steinberg [GBSt], [GB], [St] a consequence of which is that the maximal rank of a free homomorphic image of Y is 1.♠

**Corollary 11.** *Every quasi-locally free group contains a maximal Abelian subgroup elementarily equivalent to the infinite cyclic group.*

Proof: Every quasi-locally free group is non-Abelian. Let a and b be a pair of non-commuting elements. We may take g to be the commutator [a,b] in Proposition 10. (It is well known that a nontrivial commutator in a free group is neither a proper power nor a primitive element. See

p. 52 and Proposition 5.11, p. 107 of [LSc], respectively.)♠

We close this section with a proof that finitely presented, almost locally free groups have solvable conjugacy problem. To that end we require a definition.

**Definition 12.** *A group G is conjugacy separable provided to every pair $\bar{g} = (g_1, g_2)$ of elements of G which are not conjugate in G, there is a finite group $H_{\bar{g}}$ and an epimorphism $\varphi_{\bar{g}} : G \to H_{\bar{g}}$ such that $\varphi_{\bar{g}}(g_1)$ and $\varphi_{\bar{g}}(g_2)$ are not conjugate in $H_{\bar{g}}$.*

It is easy to see that conjugacy separability is a sufficient condition for a finitely presented group to have solvable conjugacy problem. It is known (see Proposition 4.8, p.26 of [LSc]) that free groups are conjugacy separable.

**Theorem 13.** *Every finitely presented almost locally free group has solvable conjugacy problem.*

Proof: Let the almost locally free group G have the finite presentation

$$<a_1, \ldots, a_m; R_1(a_1,\ldots,a_m) = \cdots = R_n(a_1,\ldots,a_m) = 1>.$$

Let $\bar{x} = (x_1,\ldots,x_m)$ be an ordered m-tuple of distinct variables. Let $u = U(a_1,\ldots,a_m)$ and $w = W(a_1,\ldots,a_m)$ be a pair of non-conjugate elements of G.

Then the $\Sigma_2$-sentence

$$\sigma: \exists \bar{x} \forall y (\bigwedge_{i=1}^{n} (R_i(\bar{x}) = 1) \wedge (W(\bar{x}) \neq y^{-1}U(\bar{x})y))$$

of $L_o$ is true in G, Could the $\Pi_2$-sentence

$$\tau: \forall \bar{x} \exists y (\bigwedge_{i=1}^{n} (R_i(\bar{x}) = 1) \wedge (W(\bar{x}) = y^{-1}U(\bar{x})y))$$

of $L_o$, logically equivalent to the negation of $\sigma$, be true in the non-Abelian free groups? If so then $\tau$ would be true in the almost locally free group G. But that contradicts the fact that $\sigma$ is true in G. Hence, $\tau$ must be false and $\sigma$ true in every non-Abelian free group. In particular, if $F_2 = <c_1, c_2; >$ is free of rank 2, then $\sigma$ is true in $F_2$. Hence, there are elements $b_1,\ldots,b_m$ of $F_2$ such that $R_i(b_1,\ldots,b_m) = 1$ for all $i = 1,\ldots,n$ and such that $U(b_1,\ldots,b_m)$ and $W(b_1,\ldots,b_m)$ are not conjugate in $F_2$.

Let F = <$b_1,\ldots,b_m$>. By Nielsen-Schreier, F is free. Now we can define an epimorphism $\psi$ : G→F by $a_i \xmapsto{\psi} b_i$, i = 1,...,m, since the relations are preserved. $\psi(u)$ and $\psi(w)$ are mapped, respectively, to the elements U($b_1,\ldots,b_m$) and W($b_1,\ldots,b_m$) which are not conjugate even in the overgroup $F_2 \geq F$. Since F is free, it is conjugacy separable so there is a finite group H and an epimorphism $\varphi$ : F → H such that $\varphi\psi(u)$ and $\varphi\psi(w)$ are not conjugate in H.

$$\begin{array}{ccc} G & \xrightarrow{\psi} & F \\ & \searrow_{\varphi\psi} & \downarrow \varphi \\ & & H \end{array}$$

By taking the composite map, we have an epimorphism $\varphi\psi$ from G onto the finite group H such that $\varphi\psi(u)$ and $\varphi\psi(w)$ are not conjugate in H. Since u and w were arbitrary, the theorem is established.♠

3. THE ABELIANIZATION OF ALMOST LOCALLY FREE GROUPS AND THE GENUS QUESTION

Observe that the Abelianization of the quasi-locally free group

$$K = <a_1, a_2, a_3, a_4; a_1^2 a_2^2 a_3^2 a_4^2 = 1>$$

contains nontrivial torsion. Indeed modulo the derived group [K,K], the element $a_1 a_2 a_3 a_4$ has order 2. We do not know whether or not the Abelianization of an almost locally free group must be torsion free. To that end the following definition has appeared earlier in the literature. (See e.g. Edmunds and Rosenberger [ER]).

**Definition 14.** *Let F be a non-Abelian free group and let [F,F] be its derived group. Define a function genus :[F,F]→ $\varpi$ by*

$$\text{genus}(z) = \begin{cases} 0 & if\ z=1 \\ \min\{N:z=[x_1,y_1]\cdots[x_N,y_N]\} & otherwise \end{cases}$$

That is, **genus** of z is the minimal number of factors necessary to express z as a product of commutators.

**Proposition 15.** *(Duncan and Howie [DH]): If g and n are positive integers and z ≠ 1 lies in the derived group of a free group F and $z^n = [x_1,y_1]\cdots[x_g,y_g]$ in F, then n < 2g.*

This result generalizes earlier partial results including the classical theorem (used in Corollary 4.1) that no nontrivial commutator can be a proper power in a free group. Suppose for each ordered pair of positive integers, we define

$$f(g,n)= \begin{cases} \infty & if \ \{\text{genus}(z):\text{genus}(z^n)\leq g\} \ \text{is unbounded} \\ \max\{\text{genus}(z): genus(z^n)\leq g\} & otherwise \end{cases}$$

Here z varies over the commutator subgroup of $F_2 = <a_1, a_2; >$; however, any non-Abelian free group will do. This is so since if $\overline{x} = (x_1,\ldots,x_g)$, $\overline{y} = (y_1,\ldots,y_g)$, $\overline{u} = (u_1,\ldots,u_k)$ and $\overline{w} = (w_1,\ldots,w_k)$, we have the $\Pi_2$-sentence $\sigma(g,n,k)$:

$$\forall \overline{x} \forall \overline{y} \exists \overline{u} \exists \overline{w} ((z^n = [x_1,y_1]\cdots[x_g,y_g]) \rightarrow (z=[u_1,w_1]\cdots[u_k,w_k])) \qquad (5)$$

either is true in every non-Abelian free group or false in every non-Abelian free group.

**Question 3:** Is $f(g,n) < \infty$ for every ordered pair (g,n) of positive integers?

Question 3 generalizes the special case g=n=2 in Edmunds and Rosenberger [ER]. By the arguments posed in [GS], $f(2,2) < \infty$ would be sufficient to rule out $K = <a_1, a_2, a_3, a_4; a_1^2 a_2^2 a_3^2 a_4^2 = 1>$ being almost locally free. Indeed, we have

**Theorem 16.** *The following three statements are pairwise equivalent.*
*(1.) $f(g,n) < \infty$ for every ordered pair (g,n) of positive integers.*
*(2.) The unrestricted direct power $F_2^\varpi$ has torsion free Abelianization.*
*(3.) Every almost locally free group has torsion free Abelianization.*

**Proof:** $(1.) \Longrightarrow (2.)$

Suppose there are elements $x_i = (x_i(k))_{k<\varpi}$, $y_i = (y_i(k))_{k<\varpi}$, $i = 1,\ldots,g$ and $z = (z(k))_{k<\varpi}$ of $F_2^\varpi$ such that $z^n = [x_1, y_1]\cdots[x_g, y_g]$. Then for all $k < \varpi$, $z(k)^n = [x_1(k), y_1(k)]\cdots[x_g(k), y_g(k)]$. So there are elements $u_1(k),\ldots,u_{f(g,n)}(k)$ and $w_1(k),\ldots,w_{f(g,n)}(k)$ such that

$$z(k) = [u_1(k), w_1(k)]\cdots[u_{f(g,n)}(k), w_{f(g,n)}(k)]$$

$\cdots$ for all $k < \varpi$. (If for some k genus(z(k)) < f(g,n), then we may take $u_j(k)$ and $w_j(k)$ to be trivial for "sufficiently large" j with $1 \leq j \leq f(g,n)$.) And so

$$z = [u_1, w_1]\cdots[u_{f(g,n)}, w_{f(g,n)}] \in [F_2^\varpi, F_2^\varpi].$$

$(2.) \Longrightarrow (3.)$

Suppose to deduce a contradiction that the Abelianization of the almost locally free group G contains nontrivial torsion. Say n > 1 is an integer and $z^n \in G' = [G,G]$ but $z \notin G'$. Suppose that $z^n = [x_1, y_1]\cdots[x_g, y_g]$. Then for all integers $k \geq 1$ the $\Pi_2$-sentence $\sigma(g,n,k)$ (see (5) must be false in G. Thus $\sigma(g,n,k)$ must be false in $F_2$.

Thus for all k $< \varpi$, we can find elements $x_i(k)$ and $y_i(k)$, i =1,...,g and z(k) of $F_2$ such that $z(k)^n = [x_1(k), y_1(k)] \cdots [x_g(k), y_g(k)]$ and such that genus(z(k)) >k. Then in $F_2^\varpi$, letting $\xi_i = (x_i(k))_{k<\varpi}$ and $\eta_i = (y_i(k))_{k<\varpi}$, i = 1,...,g, and $\varsigma = (z(k))_{k<\varpi}$, we have

$$\varsigma^n = [\xi_1, \eta_1] \cdots [\xi_g, \eta_g] \in [F_2^\varpi, F_2^\varpi]$$

but $\varsigma \notin [F_2^\varpi, F_2^\varpi]$ otherwise there would be a fixed $m < \varpi$ and elements $u_1, \ldots, u_m$ and $w_1, \ldots, w_m$ such that $\varsigma = [u_1, w_1] \cdots [u_m, w_m]$. But then $z(k) = [u_1(k), w_1(k)] \cdots [u_m(k), w_m(k)]$ for all $k < \varpi$ and genus(z(k))$\leq$ m for all $k < \varpi$. This contradiction shows that $F_2^\varpi$ having torsion free Abelianization implies that every almost locally free group has torsion free Abelianization.

(3.) $\implies$ (1.)

Suppose to deduce a contradiction that $f(g, n) = \infty$ for some ordered pair (g,n) of positive integers. Then for every positive integer k there would be elements $x_i(k)$ and $y_i(k)$, i = 1,...,g, and z(k) of $F_2$ such that $z(k)^n = [x_1(k), y_1(k)] \cdots [x_g(k), y_g(k)]$ but genus(z(k)) > k. Let D be a non-principal ultrafilter on the first infinite cardinal $\varpi$. Form the ultrapower*$F_2 = F_2^\varpi/D$. Then *$F_2$ is almost locally free (even a model of $\Sigma$). We have the elements $\xi_i = (x_i(k))_{k<\varpi}/D, \eta_i = (y_i(k)))_{k<\varpi}/D, i = 1, \ldots, g$ and $\varsigma = (z(k))_{k<\varpi}/D$ of *$F_2$ satisfying

$$\varsigma^n = [\xi_1, \eta_1] \cdots [\xi_g, \eta_g].$$

Thus $\zeta^n \in [{}^*F_2, {}^*F_2]$ and since *$F_2$ has torsion free Abelianization, we conclude $\varsigma \in [{}^*F_2, {}^*F_2]$. Hence there is a fixed m $< \varpi$ and elements $\gamma_1, \ldots, \gamma_m$ and $\delta_1, \ldots, \delta_m$, $\gamma_j = (c_j(k))_{k<\varpi}/D$ and $\delta_j = (d_j(k))_{k<\varpi}/D$, j =1,...,m such that $\zeta = [\gamma_1, \delta_1] \cdots [\gamma_m, \delta_m]$. It follows that

$$\{k \in \varpi : genus(z(k)) \leq m\} \in D$$

This contradicts the fact that for any fixed m the set

$$\{k \in \varpi : genus(z(k)) \leq m\}$$

is finite. (See Lemma 3, p.108 of [BeSl].) This contradiction shows that f(g,n) = $\infty$ is impossible.♠

**Remark 1.** *We have at least that every locally free group has torsion free Abelianization. This is so since if G is locally free and $x_i, y_i \in G$, i = 1,...,g and $z \in G$ satisfy $z^n = [x_1, y_1] \cdots [x_g, y_g]$, then the equation takes place in the finitely generated subgroup*
*$G_o = <x_1, \ldots, x_g, y_1, \ldots, y_g, z>$ and $G_o$ being free has torsion free Abelianization. So $z \in [G_o, G_o] \subseteq [G, G]$.*

We end with some observations. If $\Phi$ is the set of all universal sentences and all existential sentences of $L_o$ true in every non-Abelian free group, then the model class $\mathcal{M}(\Phi)$, i.e., the class of $\exists$-free groups, is an inductive model class since it has a set of $\Pi_2$-axioms. Thus, existentially closed objects relative to $\mathcal{M}(\Phi)$ exist in $\mathcal{M}(\Phi)$. We call such objects $\mathcal{M}(\Phi)$ e.c. groups. It is not difficult to convince oneself that every $\mathcal{M}(\Phi)$ e.c. group is almost locally free. It would be of interest to us to determine what the $\mathcal{M}(\Phi)$ e.c. groups have to "look like." For example, what - if anything - can one say about the maximal Abelian subgroups in an $\mathcal{M}(\Phi)$ e.c. group?

Our final observation is that a $\exists$-free group G has its maximal Abelian subgroups elementarily equivalent to the infinite cyclic group if and only if G is a quasi-locally free group satisfying for each prime p the following $\Pi_3$-sentence $\pi_p$ of $L_o$:

$$\forall x \exists y \forall z ((x \neq 1) \to ((xy = yx) \wedge (y \neq z^p)))$$

asserting that to every nontrivial element there is an element which commutes with it and is not a p-th power. The class of such groups is closed under the formation of non-Abelian, centralizer closed subgroups and free products but not under direct unions.